\documentclass[11pt]{amsart}
\usepackage[T1]{fontenc}
\usepackage[utf8]{inputenc}
\usepackage[usenames]{color}
\usepackage{fullpage}
\usepackage{amscd}
\usepackage{amssymb, latexsym}
\usepackage[all,2cell,ps]{xy}
\usepackage{mathdots}

\theoremstyle{plain}
\newtheorem{thm}{Theorem}[section]
\newtheorem{theorem}[thm]{Theorem}

\newtheorem{cor}[thm]{Corollary}
\newtheorem{corollary}[thm]{Corollary}

\newtheorem{lemma}[thm]{Lemma}

\newtheorem{prop}[thm]{Proposition}
\newtheorem{proposition}[thm]{Proposition}
\newtheorem{ques}[thm]{Question}

\theoremstyle{definition}
\newtheorem{de}[thm]{Definition}
\newtheorem{exm}[thm]{Example}
\newtheorem{rem}[thm]{Remark}
\newtheorem{remark}[thm]{Remark}
\newtheorem{convention}[thm]{Convention}

\newcommand{\Z}{\mathbb{Z}}
\newcommand{\N}{\mathbb{N}}

\newcommand{\id}{\mathrm{id}}

\newcommand{\Ker}{\mathop{\mathrm{Ker}}}

\newcommand{\Aut}{\mathop{\mathrm{Aut}}}

\newcommand{\Ret}{\mathop{\mathrm{Ret}}}
\newcommand{\Soc}{\mathop{\mathrm{Soc}}}
\newcommand{\Ann}{\mathop{\mathrm{Ann}}}
\newcommand{\nil}{\mathop{\mathrm{nil}}}

\numberwithin{equation}{section}

\begin{document}
\title{Cocyclic braces and indecomposable cocyclic solutions of the Yang-Baxter equation}

\author{P\v remysl Jedli\v cka}
\author{Agata Pilitowska}
\author{Anna Zamojska-Dzienio}

\address{(P.J.) Department of Mathematics, Faculty of Engineering, Czech University of Life Sciences, Kam\'yck\'a 129, 16521 Praha 6, Czech Republic}
\address{(A.P., A.Z.) Faculty of Mathematics and Information Science, Warsaw University of Technology, Koszykowa 75, 00-662 Warsaw, Poland}

\email{(P.J.) jedlickap@tf.czu.cz}
\email{(A.P.) A.Pilitowska@mini.pw.edu.pl}
\email{(A.Z.) A.Zamojska@mini.pw.edu.pl}

.

\date{\today}

\begin{abstract}
We study indecomposable involutive set-theoretic solutions of the Yang-Baxter equation with cyclic permutation groups (cocyclic solutions). In particular, we show that there is no one-to-one correspondence between indecomposable cocyclic solutions and cocyclic braces which contradicts recent results in \cite{Rump21}. 
\end{abstract}

\maketitle
\section{Introduction}
The Yang-Baxter equation is a fundamental equation occurring in integrable models in statistical mechanics and quantum field theory~\cite{Jimbo}.
Let $V$ be a vector space. A {\em solution of the Yang--Baxter equation} is a linear mapping $r:V\otimes V\to V\otimes V$ such that
\[
(id\otimes r) (r\otimes id) (id\otimes r)=(r\otimes id) (id\otimes r) (r\otimes id).
\]
Description of all possible solutions seems to be extremely difficult and therefore
there were some simplifications introduced (see e.g. \cite{Dr90}).

Let $X$ be a basis of the space $V$ and let $\sigma:X^2\to X$ and $\tau: X^2\to X$ be two mappings.
We say that $(X,\sigma,\tau)$ is a {\em set-theoretic solution of the Yang--Baxter equation} if
the mapping $x\otimes y \mapsto \sigma(x,y)\otimes \tau(x,y)$ extends to a solution of the Yang--Baxter
equation. It means that $r\colon X^2\to X^2$, where $r=(\sigma,\tau)$ satisfies the \emph{braid relation}:

\begin{equation}\label{eq:braid}
(id\times r)(r\times id)(id\times r)=(r\times id)(id\times r)(r\times id).
\end{equation}

A solution is called {\em non-degenerate} if the mappings $\sigma_x:=\sigma(x,\_)$ and $\tau_y:=\tau(\_\,,y)$ are bijections,
for all $x,y\in X$.
A solution $(X,\sigma,\tau)$ is {\em involutive} if $r^2=\mathrm{id}_{X^2}$. In the involutive case, the operation $\tau$ can be expressed by means of the operation $\sigma$ (see \cite[Proposition 1]{Rump05}).

\begin{convention}
All solutions, we study in this paper, are set-theoretic, non-degenerate and involutive so we will call them simply \emph{solutions}.
\end{convention}

\emph{The permutation group} $\mathcal{G}(X)=\left\langle \sigma_x\colon x\in X\right\rangle$ of a solution $(X,\sigma,\tau)$ is the subgroup of the symmetric group $S(X)$ generated by mappings $\sigma_x$, with $x\in X$. The group $\mathcal{G}(X)$ is also called \emph{the involutive Yang-Baxter group} (IYB group) associated to the solution $(X,\sigma,\tau)$. A solution $(X,\sigma,\tau)$ is \emph{indecomposable} if the permutation group  $\mathcal{G}(X)$ acts transitively on $X$.

It was already observed in~\cite{ESS} that indecomposable solutions
form an important class of all solutions. Since then, many articles on indecomposable solutions appeared, for instance (e.g. \cite[Problem 4]{Vendr}, \cite{CCP,CPR20,JPZ20,Rump}).

Solutions~$(X,\sigma,\tau)$ with $\mathcal{G}(X)$ cyclic are called {\em cocyclic}. A classification of such solutions was recently given by Rump~\cite{Rump21}. In this paper we focus on indecomposable cocyclic solutions. Among the others, we
reformulate the main result~\cite[Theorem 1]{Rump21} here as Theorem~\ref{thm:RumpMain}. 
According to Rump's classification, a complete system of invariants
are: the size of the solution and the index of a subgroup of~$\mathcal{G}(X)$. 
In particular, this classification claims that there exist, for each prime~$p$,
exactly two non-isomorphic indecomposable solutions with their permutation groups isomorphic to~$\Z_{p^2}$.

However, there is an earlier classification by Castelli, Pinto and Rump~\cite{CPR20} of all the indecomposable solutions of size~$p^2$ with abelian permutation groups.
For each prime~$p$, there exist $p+1$ such solutions up to isomorphism, among which one has
its permutation group isomorphic to~$\Z_p^2$ and the other~$p$ solutions
have their permutation groups isomorphic to~$\Z_{p^2}$.
This result has been confirmed by the authors in~\cite{JPZ20} using a different technique.

The main motivation for this paper was therefore to investigate the discrepancy between those
two results cited above. We analyzed several results about indecomposable solutions from the point of view of their permutation groups and this has lead us to a final conclusion: our results are summarized in Theorem~\ref{th:main}.
For an indecomposable cocyclic solution, to constitute a complete system of invariants one needs three parameters not two and, as a consequence, there is no one-to-one correspondence between cocyclic braces and indecomposable cocyclic solutions. This shows that Theorem 1 and its Corollary in \cite{Rump21} are incorrect.

The paper is organized as follows: in Section~2 we give basic definitions used
in contemporary studies of involutive set-theoretic solutions of Yang-Baxter equation. We also investigate all the
9-element indecomposable solutions, as the smallest example where the two theories mentioned above differ. This investigation will eventually bring us to a trail of a mistake in \cite{Rump21}. In Section~3 we present a construction by Bachiller, Cedó and Jespers~\cite{BCJ16} that allows one to construct all the solutions with a given permutation group. 
We first analyze prime power size solutions and then we use Chinese remainder theorem to
capture all finite sizes.
At the end, in Section~4, we use this construction to enumerate all the indecomposable solutions with their permutation group cyclic.

\section{Preliminaries}

In this section we give all the necessary notions that we need to analyse both papers~\cite{Rump21} and~\cite{CPR20}. We give also an example where the results of both the papers differ.

Rump in \cite[Definition 2]{Rump07A} introduced the notion of a brace and showed the correspondence between such structures and solutions. Here we use an equivalent definition formulated by Ced\'{o}, Jespers and Okni\'{n}ski.

\begin{de}\cite[Definition 2.2]{CJO14}
 An algebra $(B,+,\circ)$ is called a {\em left brace} if $(B,+)$ is
 an abelian group, $(B,\circ)$ is a group and the operations satisfy, for all $a,b,c\in B$,
 \begin{equation}\label{lb}
  a\circ b+a\circ c = a\circ(b+c)+a.
 \end{equation}
 \end{de}
 \noindent
The inverse element to~$a$ in the group $(B,\circ)$ shall be denoted by $\bar a$. Notice that identity elements in both groups $(B,+)$ and $(B,\circ)$ coincide so we will indicate them by $0$.

A left brace $(B,+,\circ)$ is said to be \emph{abelian} if its adjoint (or \emph{multiplicative}) group $(B,\circ)$ is abelian, \emph{cyclic} if the additive group $(B,+)$ is cyclic and it is called \emph{cocyclic} if its adjoint group $(B,\circ)$ is cyclic. In case that both groups are cyclic, the brace is then called \emph{bicyclic}.

For a left brace $(B,+,\circ)$ there is an action $\lambda\colon(B,\circ)\rightarrow\Aut{(B,+)}$, where $\lambda(a)=\lambda_a$ and for all $b\in B$, $\lambda_a(b)=a\circ b-a$. The kernel of the group homomorphism $\lambda$ is called the \emph{socle} of the left brace $(B,+,\circ)$ and denoted by $\Soc(B)$, i.e.
\begin{equation*}
\Soc(B)=\Ker(\lambda)=\{a\in B\colon \lambda_a=\id\}=\{a\in B\colon a\circ b=a+b \text{ for all } b\in B\}.
\end{equation*}
$\Soc(B)$ is a normal subgroup of the group $(B,\circ)$ and is an ideal of the left brace. Hence the quotient $(B,\circ)/\Soc(B)$ of the multiplicative group is also the quotient of the additive group $(B,+)$ and the factor brace $B/\Soc(B):=(B,+,\circ)/\Soc(B)$ by the ideal $\Soc(B)$ is well defined.

A right brace is defined similarly, replacing Property \eqref{lb} by $b\circ a+c\circ a = (b+c)\circ a+a$. If $(B,+,\circ)$ is both a left and a right brace then one says that it is a (\emph{two-sided}) brace. A brace is called {\em trivial} if $a+b=a\circ b$, for all $a,b\in B$.

For a left brace $(B,+,\circ)$ one can define another operation $\ast$ as follows
\begin{equation*}
a\ast b=a\circ b-a-b,
\end{equation*}
for $a,b\in B$.
A brace $(B,+,\circ)$ is two-sided if and only if $(B,+,\ast)$ is a radical ring \cite[Proposition 1]{CJO14}. In particular, every abelian left brace is two-sided.

The \emph{associated solution} $(B,\sigma,\tau)$ to a left brace $(B,+,\circ)$ is defined on a set $B$ as follows: $\sigma_x(y):=\lambda_x(y)$ and $\tau_y(x):=\lambda^{-1}_{\lambda_x(y)}(x)$ for $x,y\in B$.

 %There is one-to-one correspondence between factor braces by $\Soc(B)$ and retractions of solutions associated with the braces.
 
In \cite[Section 3.2]{ESS} Etingof, Schedler and Soloviev introduced, for each solution $(X,\sigma,\tau)$, the equivalence relation $\sim$ on the set $X$: for each $x,y\in X$
\[
x\sim y\quad \Leftrightarrow\quad \sigma_x=\sigma_y.
\]
They showed that the quotient set $X/\mathord{\sim}$ can be again endowed
with a structure of a solution. They call such a solution the {\em retraction} of the solution~$(X,\sigma,\tau)$ and denote it by
$\Ret(X)$. One can also define \emph{iterated retraction} in the following way: $\Ret^0(X,\sigma,\tau):=(X,\sigma,\tau)$ and
$\Ret^k(X,\sigma,\tau):=\Ret(\Ret^{k-1}(X,\sigma,\tau))$, for any natural number $k>1$.
A solution $(X,\sigma,\tau)$ is called a \emph{multipermutation solution of level $m$} if $m$ is the least nonnegative integer such that
$|\Ret^{m}(X,\sigma,\tau)|=1.$ In such case, we will also say that a solution is \emph{of a multipermutational level $m$}. %A solution is \emph{retractable} if such an $m$ exists.

\begin{proposition}\label{prop:solsoc}\cite[Proposition 7]{Rump07A} $($see also \cite[Lemma 3]{CJO14}$)$\label{prop:ret}
Let $(B,+,\circ)$ be a nonzero left brace. The solution associated to the factor left brace $B/\Soc(B)$ is equal to the retraction of the solution associated to $(B,+,\circ)$.
\end{proposition}

There are many different methods to construct braces. Etingof, Schedler and Soloviev in \cite[Appendix]{ESS} introduced the notion of a $T$-\emph{structure} on an abelian group $A$. By \cite[Theorem A7]{ESS} one can deduce that cyclic braces correspond to $T$-structures on cyclic groups. Rump in \cite{Rump07B} used this result to obtain a wide family of cyclic braces constructed on rings with an identity and with a cyclic additive group. See \cite[Proposition 1]{Rump07B} and \cite[Theorem 1]{Rump07B} for the particular case of the ring $(\Z_n,+,\cdot)=\mathbb{Z}/(n)$.

By results of \cite[Example 3, Appendix]{ESS} Rump's construction can be easily generalized. 
\\
The ring multiplication in any ring $(R,+,\cdot)$ will be consequently indicated by juxtaposition.
\begin{proposition}\label{prop:brace}
Let $(R,+,\cdot)$ be a ring with identity $1$,  $J(R)$ be the Jacobson radical of the ring and $r\in J(R)$. Let for $a,b\in R$, 
\begin{align}
a\circ b&:=a+b+arb. 
\end{align}

Then $(R,+,\circ)$ is a two-sided brace with  $\Soc(R)=\Ann(r)$.
  \end{proposition}
\begin{proof}

Directly by Etingof, Schedler and Soloviev \cite{ESS},   $(R,\circ)$
is a group with $0$ as the unit and $\bar a:=-a(1+r a)^{-1}$ as the inverse of an element $a$. 

Straightforward calculations show that  $(R,+,\circ)$ is a two-sided brace:
\[a\circ(b+c)+a=a+b+c+ar(b+c)+a
 = a+b+arb+a+c+arc=a\circ b+a\circ c,
\]
and
\[(b+c)\circ a+a=b+c+a+(b+c)ra+a
 = b+a+bra+c+a+cra=b\circ a+c\circ a.
\]
Now, for all $b\in R$,
\[ a\in\Soc(R) \Leftrightarrow
 \lambda_a(b)=b \Leftrightarrow a+b=a\circ b
 \Leftrightarrow arb=0 \Leftrightarrow
 a\in \Ann(r). \qedhere
\]
\end{proof}

Note that the construction from Proposition \ref{prop:brace} is valid for any nilpotent element $r\in R$ (see e.g. \cite[Lemma 4.4]{KN07}). 

\begin{de}\cite{CGS17}
Let $(B,+,\circ)$ be a left brace. It is \emph{left nilpotent} if $B^n = 0$ for some positive integer $n$, where $B^{n+1} = B\ast B^n$ and it is \emph{right nilpotent} if $B^{(m)} = 0$ for some $m\in\mathbb{Z}^+$, where $B^{(1)}=B$ and $B^{(m+1)} = B^{(m)}\ast B$. 
\end{de}
The least $n\in\Z^+$ such that $B^{n+1} = 0$ is called \emph{the left nilpotency degree} of a left brace $(B,+,\circ)$ or one shortly says that 
$(B,+,\circ)$ is a \emph{left $n$-nilpotent} left brace. Similarly, if $m\in\Z^+$ is the least number such that $B^{(m+1)} = 0$, then $(B,+,\circ)$ is a \emph{right $m$-nilpotent} left brace or a left brace of the right nilpotency degree equal to $m$.
There is a direct connection between the right nilpotency degree of a left brace and the multipermutational level of its associated solution.

\begin{proposition}\cite[Proposition 6]{CGS17}\label{prop:CGS}
Let $(B,+,\circ)$ be a nonzero left brace and let $(B,\sigma,\tau)$ be its associated solution. Then $(B,\sigma,\tau)$ is a solution of a multipermutational level $m$ if and only if $(B,+,\circ)$ is right $m$-nilpotent.
\end{proposition}

If a left brace is abelian, then the left nilpotency and the right nilpotency
are the same notions and we speak simply of {\em nilpotency}.
In the finite case, every abelian left brace is nilpotent (\cite[Proposition 3]{CJO14}).

Let $(R,+,\cdot)$ be a commutative ring with identity $1$ and $J(R)$ the Jacobson radical of the ring. Let $(R,+,\circ)$ be the two-sided brace defined in Proposition \ref{prop:brace} for some $r\in J(R)$. 
 Then, for $a,b\in R$ we have $a\ast b=rab$. Since the operation $\ast$ is commutative and associative this implies, for arbitrary $k\in \mathbb{N}$
 $$R^k=R^{(k)}=r^{k-1} R.$$
We then immediately obtain:
 
\begin{proposition}\label{prop:ring_1}
Let $(R,+,\cdot)$ be a commutative ring with identity $1$ and $J(R)$ the Jacobson radical of the ring. Let $(R,+,\circ)$ be the brace defined in Proposition \ref{prop:brace}, for some $r\in J(R)$. 
 Then, 
 \begin{enumerate}
  \item The brace~$(R,+,\circ)$ is nilpotent if and only if the element $r$ is nilpotent.\\
   In such a case the nilpotency degrees of $r$ and $(R,+,\circ)$ are equal. %$m$ such that $R^{(m+1)}=0$ and $R^{(m)}\neq 0$.
   \item The associated solution to the brace $(R,+,\circ)$  is a multipermutation solution if and only if the element $r$ is nilpotent. The multipermutational level is equal to degree of nilpotency of $r$. 
 \end{enumerate}
\end{proposition}
\noindent

Our work focuses on finite cocyclic solutions that are closely connected to the following example of cyclic braces.

\begin{exm}\label{exm:coc_br} (See also \cite[4.1]{KN08}.)
Let $n=p_1^{k_1}\cdots p_s^{k_s}$ for some prime pairwise distinct numbers $p_1,\ldots,p_s$ and positive integers $k_1, \ldots, k_s$. 
The nilradical $\nil(\Z_n)$ of the ring $(\Z_n,+,\cdot)$ consists precisely of the elements of the form $rp_1^{t_1}\cdots p_s^{t_s}$, for $1\leq t_1\leq k_1, \ldots, 1\leq t_s\leq k_s$ and $r\in\Z_n^*$.

Let $t\in\mathbb{Z}$ be such that $t \bmod{n} \in \nil(\Z_n)$. By Proposition \ref{prop:brace}, $(\Z_n,+,\circ)$ with
$$a\circ b=a+ b+t ab$$
is a two-sided brace with $\Soc(\Z_n)=\frac{n}{t}\Z_n\cong \Z_t$. 

By Proposition \ref{prop:ring_1}, the associated solution to such brace is a solution of a multipermutational level equal to $\max(\lceil\dfrac{k_1}{t_1} \rceil,\lceil\dfrac{k_2}{t_2} \rceil,\dots,\lceil\dfrac{k_s}{t_s} \rceil)$. 

We will denote the brace $(\Z_n,+,\circ)$, obtained as described above, by $B_t(n)$.

\noindent
Note that braces constructed for $t=p_1^{t_1}\cdots p_s^{t_s}$ and $rt$, for $r\in\Z_n^*$, are isomorphic via the mapping $f(x)=r^{-1}x$.
This justifies the introduction of the notation: 
$$
\mathcal{N}(\Z_n):=\{p_1^{t_1}\cdots p_s^{t_s}\colon t_i\in\{1,\ldots k_i\} \text{ for } i\in\{1,\ldots s\}\},
$$
which we use in the following sections.
%In case $t\in\mathcal{N}(\Z_n)$ we will denote the brace $(\Z_n,+,\circ)$, obtained as described above, by $B_t(n)$.
\end{exm}

Braces of the form $B_t(n)$ are actually the only finite abelian cyclic braces \cite[Theorem 1]{Rump07B}. By \cite[Proposition 6]{Rump07B} this class contains almost all finite cocyclic braces.

For each solution $(X,\sigma,\tau)$, Etingof, Schedler and Soloviev \cite[Section 2]{ESS} introduced its \emph{structure group} $(G(X),\cdot)$ generated by elements of $X$ with defining relations $x\cdot y=\sigma_x(y)\cdot \tau_y(x)$. 
$G(X)$ has also a natural structure of a left brace with an additive group being a free abelian group with basis $X$ (see e.g \cite{CJO14,CGS18}). This left brace is denoted by $A_X$.
Furthemore, let $A(X)$ be the factor brace $A_X/\mathrm{Soc}(A_X)$.

%Furthemore, by $A(X)$ we denote the factorbrace $A_X/\mathrm{Soc}(A_X)$.

%??? Nilpotency degree of $A(X)$ = multipermutation level of~$X$ ???

Cocyclic solutions have been recently studied in~\cite{Rump21}. By construction,
the adjoint group of~$A(X)$ is isomorphic to the permutation group~$\mathcal{G}(X)$. Then, knowing the permutation group of a solution and using an algorithm described in \cite[Section 1]{Rump21}, one is able to reconstruct the solution~$(X,\sigma,\tau)$ from the brace $A(X)$. This leads to the following classification:

\begin{theorem}\cite[Theorem 1]{Rump21}\label{thm:RumpMain}
For the isomorphism class of an indecomposable cocyclic solution~$(X,\sigma,\tau)$, the order $n$ and the socle index $s=|A(X)/\mathrm{Soc}(A(X))|$ form a complete system of invariants. For a given order~$n$, the possible socle indices are the
divisors~$s$ of~$n$ which satisfy $p\mid n \Rightarrow ps\mid n$, for all primes $p$ and $8\mid n \Rightarrow 4s\mid n$.
\end{theorem}

An immediate consequence of the theorem is

\begin{cor}\cite{Rump21}\label{cor:Rump}
Up to isomorphism, there is a one-to-one correspondence between finite indecomposable cocyclic solutions and finite cocyclic braces.
\end{cor}

This corollary was quite unexpected, as Rump himself admits in the introduction of~\cite{Rump21}.
Nevertheless, as we already announced in the introduction of this paper, both of the results are incorrect.
We illustrate it on the smallest counterexample.

\begin{exm}\label{exm:sol}
Consider the ring~$(\Z_9,+,\cdot)$. According to~\cite[Corollary 22]{CPR20} and \cite[Corollary 4.8]{JPZ20}, there exist three indecomposable $9$-element solutions $(X,\sigma,\tau)$ 
with their permutation group cyclic, namely for $r\in\{0,3,6\}$, the set $X=\{0,1,\ldots,8\}$ with
\[\sigma_{x}(y)=rx+y+1 \qquad \text{and}\qquad \tau_y(x)=x-1-r(y+1).\]
They are not isomorphic since, for each~$x\in X$, $\sigma_{\sigma_x(x)}=\sigma_x^{r+1}$.

It can be computed that the brace $A(X)$ is then isomorphic to the brace
$B_{9-r}(9)$ described in Example \ref{exm:coc_br}. Of course, for $r=0$, it is a trivial brace $B_9(9)$, whereas for $r\in\{3,6\}$ one obtains the brace $B_3(9)\cong B_6(9)$ (see also Example \ref{exm:n=9}).
\end{exm}

This example indicates (a thorough proof comes in the next section though)
that there is a mistake in the proof of Theorem~\ref{thm:RumpMain}. Indeed, there
is one, a~very subtle one: on page~473 the autor correctly supposes that
the multiplicative group can be identified with~$\Z_n$. But then the proof of
the main theorem on page~475 starts with: ``\emph{Up to isomorphism, the solution $X$
is determined by the cocyclic brace $A(X)$}'' (which is unique up to isomorphism). However,
we cannot use the ``\emph{up to isomorphism}'' argument anymore since
we have already chosen an isomorphism by identifying the elements of~$A(X)$ with
the elements of the group~$\Z_n$.

The rest of the proof is correct. Therefore we can say that a cocyclic solution
is uniquely determined by the brace $A(X)$ and by the choice of the isomorphism
$(A(X),\circ)\cong \Z_n$, that means by the choice of a generator of~$(A(X),\circ)$. Nevertheless, different choices of the generators may still yield isomorphic solutions. Finding a correct answer for the isomorphism question using the tools
of~\cite{Rump21} does not seem to be straightforward and therefore, in the next section, we introduce results from yet another paper.

\section{Cocyclic braces and solutions}

In this Section our aim is to construct all finite  indecomposable solutions with cyclic permutation group. 
In \cite{BCJ16} Bachiller, Ced\'{o} and Jespers present the construction of all finite solutions $(X,\sigma,\tau)$ of the YBE such that the permutation group $\mathcal{G}(X)$ is isomorphic, as a left brace, to a given finite left brace $(B,+,\circ)$. 
We shall not repeat their construction here, as our case is much less complicated than the general one.
If the permutation group $\mathcal{G}(X)$ is abelian and acts on $X$ transitively, then it acts regularly on $X$ and therefore its order is equal to $|X|$ and, simultaneously, to $|B|$. 

By \cite[Theorem 19]{S18}, a finite left brace of cardinality $n=p_1^{k_1}\ldots p_s^{k_s}$, for some pairwise distinct prime  numbers $p_1,\ldots,p_s$ and $k_1, \ldots, k_s\in \Z^+$, with multiplicative group isomorphic to an abelian group $\Z_n$ is a direct product of left braces 
whose cardinalities are powers of prime numbers.  By \cite[Section 5]{Rump07B}, every cyclic brace of a cardinality equal to a power of some odd prime is bicyclic.  On the other hand, by \cite[Proposition 5.4]{BCJ16}, every cocyclic left brace $(B,+,\circ)$ of order $4\neq n=p^k$, for some prime $p$ and $k\geq 1$, has a cyclic additive group $\Z_{p^k}$. According to Rump's classification \cite[Theorem 1]{Rump07B} (see also \cite[Section 5]{BCJ16}), such braces are equal to $B_t(p^k)$ (up to isomorphism), for $t\in \mathcal{N}(\Z_{p^k})$ with assumption that $t\neq 2$ when $p^k=4$. 
Four-element left braces form an exception in the classification: there exist two trivial ones, one cyclic which is not cocyclic and one cocyclic which is not cyclic \cite{Rump21}. See also a clasification of commutative radical rings given in \cite[Proposition 4.1.7]{KN08}.

This implies (see \cite[Remark 5.5]{BCJ16})  that a cocyclic left brace of cardinality $n=2^{k_1}\cdot p_2^{k_2}\ldots p_s^{k_s}$ with $0\leq k_1\neq 2$ is isomorphic to the following direct product 
\[
B_{t_1}(2^{k_1})\times B_{t_2}(p_2^{k_2})\times \ldots\times B_{t_s}(p_s^{k_s}), 
\]
with some $t_i\in \mathcal{N}(\Z_{p_i^{k_i}})$, for $1\leq i\leq s$, where $t_1\neq 2$ or $k_1=1$.
%$t_1\leq k_1, t_2\leq k_2,\ldots, t_s\leq k_s$. 

Any such direct product is isomorphic to the brace 
described in Example \ref{exm:coc_br}.

\begin{lemma}\label{lem:isobrace}
Let $n=p_1^{k_1}\cdots p_s^{k_s}$ for some pairwise distinct prime numbers $p_1,\ldots,p_s$ and $k_1, \ldots, k_s\in \Z^+$. 
Let $1<t_i\leq p_i^{k_i}$ be a power of~$p_i$, for $1\leq i\leq s$, and let $t=t_1\cdots t_s$. Then the braces $B_t(n)$ and 
$B_{t_1}(p_1^{k_1})\times B_{t_2}(p_2^{k_2})\times \ldots\times B_{t_s}(p_s^{k_s})$ are isomorphic.
\end{lemma}
\begin{proof}
Indeed, let $f\colon \Z_n\to \Z_{p_1^{k_1}}\times \Z_{p_2^{k_2}}\times \ldots\times \Z_{p_s^{k_s}}$
be such that
\[\quad a\mapsto (a\frac{t}{t_1}\pmod{p_1^{k_1}},\ldots,a\frac{t}{t_i}\pmod{p_i^{k_i}},\ldots,a\frac{t}{t_s}\pmod{p_s^{k_s}}).
\]
Since for each $1\leq i\leq s$, the element $\frac{t}{t_i}$ is a generator of the group $\Z_{p_i^{k_i}}$, the mapping $f$ is a homomorphism of additive groups. Next, for $a,b\in \Z_n$
\begin{align*}
&f(a)\circ f(b)=(a\frac{t}{t_1}\pmod{p_1^{k_1}},\ldots ,a\frac{t}{t_s}\pmod{p_s^{k_s}})\circ
  (b\frac{t}{t_1}\pmod{p_1^{k_1}},\ldots ,b\frac{t}{t_s}\pmod{p_s^{k_s}})=\\
  &(a\frac{t}{t_1}\pmod{p_1^{k_1}}\circ b\frac{t}{t_1}\pmod{p_1^{k_1}},\ldots ,a\frac{t}{t_s}\pmod{p_s^{k_s}}\circ b\frac{t}{t_s}\pmod{p_s^{k_s}})=\\
&(a\frac{t}{t_1}+b\frac{t}{t_1}+t_1\cdot a\frac{t}{t_1}\cdot b\frac{t}{t_1}\pmod{p_1^{k_1}},\ldots ,a\frac{t}{t_s}+b\frac{t}{t_s}+t_s\cdot a\frac{t}{t_s}\cdot b\frac{t}{t_s}\pmod{p_s^{k_s}})=\\
&((a+b+tab)\frac{t}{t_1}\pmod{p_1^{k_1}},\ldots,(a+b+tab)\frac{t}{t_s}\pmod{p_s^{k_s}})=f(a+b+tab)=f(a\circ b).
\end{align*}
Hence,  $f$ is a homomorphism of braces, too. Finally, if
\begin{align*}
 f(a)=(a\frac{t}{t_1}\pmod{p_1^{k_1}},\ldots,a\frac{t}{t_s}\pmod{p_s^{k_s}})=
  (b\frac{t}{t_1}\pmod{p_1^{k_1}},\ldots,b\frac{t}{t_s}\pmod{p_s^{k_s}})=f(b)
\end{align*}
we have 
$$\forall(1\leq i\leq s)\quad p_i^{k_i}|(a-b)\quad \Rightarrow\quad n|(a-b)\quad \Rightarrow\quad a=b,$$
and $f$ is a bijection.
\end{proof}

As a corollary one obtains generalization of \cite[Proposition 5.4]{BCJ16} (see also \cite[Proposition 6]{Rump07B}).
\begin{corollary}\label{cor:1}
Let $n\in\Z^+$ and $t\in\mathcal{N}(\Z_n)$ be such that $4|t$ whenever $4|n$. Then the brace $B_t(n)$ is bicyclic.
\end{corollary}

\begin{remark}
The multiplicative group of the brace $B_2(2^{k+1})$, for $k\in\Z^+$, is isomorphic to the (non-cyclic) group $\Z_2\times\Z_{2^k}$. 
\end{remark}

A cocyclic left brace of cardinality $n=2^2\cdot p_2^{k_2}\ldots p_s^{k_s}$
also splits as a direct product, namely
\[
\tilde B_{t_1}(4)\times B_{t_2}(p_2^{k_2})\times \ldots\times B_{t_s}(p_s^{k_s}), \]
with some $t_i\in \mathcal{N}(\Z_{p_i^{k_i}})$, for $1\leq i\leq s$, where
$\tilde B_{t_1}(4)$ is a cocyclic brace obtained from $B_{t_1}(4)$ by
swapping the operations $+$ and $\circ$.

\bigskip

To construct solutions with the permutation group isomorphic to 
the associated group of a given left brace, it is enough to find two special mappings: $\eta\colon X\to B$ and $\varrho\colon B\to S(X)$ as described in \cite[Section 2]{BCJ16}. In the case when the left brace $(B,+,\circ)$ is cocyclic of order $p^k$, for some prime $p$ with $k\in \Z^+$, and solutions are indecomposable,  this construction can be significantly simplified.

\subsection{Indecomposable cocyclic solutions of order $p^k$}
First of all, the construction in~\cite{BCJ16}, Section~3, involves an intersection of subgroups
$\bigcap_I\bigcap_{J_i}\bigcap_b bK_{i,j}b^{-1}=\{0\}$.
Since any subgroup $K$ of the abelian group $(B_t(p^k),\circ)$ is normal we of course have $bK_{i,j}b^{-1}=K_{i,j}$, for any $b\in B$. Furthermore, to obtain the trivial intersection of a family of subgroups of the cyclic group $(B_t(p^k),\circ)$, we have to assume that at least one $K_{i,j}$ is equal to $\{0\}$.

Moreover, in \cite{BCJ16}, Theorem 3.1 starts with $X=\sqcup_I\sqcup_{J_i} B/K_{i,j}$. Given
the equality $|B|=|\mathcal{G}(X)|=|X|$
and the fact $K_{i,j}=\{0\}$, for some $i\in I$ and some $j\in J_i$,
we necessarily have $|I|=|J_{i}|=1$.
%This forces, for the cocyclic brace $(B_t(p^k),+,\circ)$,
%$$|\mathcal{G}(X)|=|X|=|B_t(p^k)|=p^k$$
%and 
Since $(B,+)$ has to be generated by $\bigcup\bigcup I$, there exists
some $a\in  \Z_{p^k}^*$ such that
$$\eta\colon \Z_{p^k}\to \Z_{p^k}; \quad \eta(b)=\lambda_b(a),$$
and 
$$\varrho\colon  \Z_{p^k}\to S( \Z_{p^k}); \quad \varrho(c)\colon \Z_{p^k}\to \Z_{p^k}; \quad \varrho(c)(b)=c\circ b.$$
Hence
$$\varrho(\eta(b))(c)=
\varrho(\lambda_b(a))(c)=\lambda_b(a)\circ c=
(b\circ a-b)\circ c=$$
$$(b+a+tab-b)+c+t(a+tab)c=a+tab+c+t(a+tab)c.$$
This implies the following reformulation of \cite[Theorem 3.1]{BCJ16}.

\begin{thm}\label{thm:constr}
Let $p$ be a prime number and $k\in\Z^+$. Let $t\in\Z$ such that
$t\;\bmod {p^k}\;\in\mathrm{nil}(\Z_{p^k})$
and $(B_t(p^k),+,\circ)$ be a cocyclic brace of order $p^k$. 

For any indecomposable solution $(X,\sigma,\tau)$ with the permutation group $\mathcal{G}(X)$ isomorphic, as left braces, to $(B_t(p^k),+,\circ)$, 
there exists $a\in \Z_{p^k}^*$, so that $(X,\sigma,\tau)$
is isomorphic to the solution defined on the set $\Z_{p^k}$ with 
$$\sigma_b(c):=a+c+tab+tac+t^2abc$$
for all $b,c\in \Z_{p^k}$.

%$a\in Orb(t)$ such that $(B_t(p^k),+)=\langle Orb(t)\rangle$.

\end{thm}
Let us denote such a solution by $\mathcal{K}(p^k,t,a)$.
\begin{rem}\label{rem:ess} \begin{enumerate}
\item For $k=1$ one obtains \cite[Theorem 2.13]{ESS}: the indecomposable solution with a prime number of elements is a cyclic permutation solution, i.e. a multipermutation solution of level $1$ with $\sigma_x$ being a cycle. There is only one (up to isomorphism) such solution for each prime. 
\item If $t\in\mathcal{N}(\Z_{p^k})$ and
 $r\in \Z^*_{p^k}$ then, similarly as for braces in Example \ref{exm:coc_br}, we obtain an isomorphic solutions 
 $\mathcal{K}(p^k,t,a)$ and $\mathcal{K}(p^k,rt,ar^{-1})$
 via the mapping $\Phi(x)=r^{-1}x$. Thus we usually restrain our focus to $t\in\mathcal{N}(\Z_{p^k})$.
\end{enumerate}
\end{rem}

\cite[Theorem 4.1]{BCJ16} characterizes when two solutions, with their permutation groups isomorphic to the same brace, are isomorphic.
Again, the assumptions $|I|=1$ and $K_{i,j}=\{0\}$ simplify the theorem substantially.
\begin{thm}
Let $a,a'\in \Z_{p^k}^*$. The solutions $\mathcal{K}(p^k,t,a)$ and $\mathcal{K}(p^k,t,a')$ are isomorphic if and only if there exist an automorphism $\Psi$ of the left brace $(B_t(p^k),+,\circ)$ and an element $z\in \Z_{p^k}$ such that 
$$\Psi(a)=\lambda_z(a').$$
\end{thm}
Each automorphism $\Psi$ of the cyclic group $\Z_{p^k}$ is of the form 
$\Psi(x)=\alpha x$,
where $\alpha \in \Z_{p^k}^*$. 
Moreover, each automorphism $\Psi$ of the left brace $(B_t(p^k),+,\circ)$ must satisfy additionally, for every $x,y\in \Z_{p^k}$, the following:
\begin{multline*}
\alpha x+\alpha y+\alpha^2 txy=\alpha x\circ \alpha y=\Psi(x)\circ \Psi(y)=\\
\Psi(x\circ y)=\alpha(x\circ y)=\alpha(x+y+txy)=
\alpha x+\alpha y+\alpha txy.
\end{multline*}
This implies that, for any $x,y\in \Z_{p^k}$, $\alpha txy\equiv\alpha^2 txy\; \pmod{p^k}$ or equivalently, $p^k\mid txy\alpha(\alpha-1).$
Hence, for $x=y=1$ and remembering that $t=p^w$ for some $w\in\left\{1,\ldots ,k\right\}$, $p^{k-w}|\alpha(\alpha-1).$
Since $\mathrm{gcd}(\alpha,p)=1$ we obtain $p^{k-w}|(\alpha-1)$ which means that $\alpha \equiv 1\; \pmod{p^{k-w}}$. Finally, we obtain that
$$\exists(h\in \Soc(B_t(p^k)))\; \alpha=1+h.$$
Summarizing, 
\begin{cor}
The solutions $\mathcal{K}(p^k,t,a)$ and $\mathcal{K}(p^k,t,a')$ are isomorphic if and only if there exist
$h\in \Soc(B_t(p^k))$ and $z\in \Z_{p^k}$ such that 
\begin{equation}\label{eq:isosol}
(1+h)a=a'(1+zt).
\end{equation}
\end{cor}
In particular, for $h=0$, $a$ and $a'$ are in the same orbit of the action $\lambda$.
%$$\lambda\colon (B_t(p^k),\circ)\to Aut(\Z_{p^k},+),$$
%$$b\mapsto \lambda_b\colon \Z_{p^k}\to \Z_{p^k}; \quad \lambda_b(a)=b\circ a-b.$$
Thus all solutions $\mathcal{K}(p^k,t,a)$ and $\mathcal{K}(p^k,t,a')$ with $a$ and $a'$ in the same orbit of the action $\lambda$ are isomorphic.
\begin{cor}\label{cor:iso}
The solutions $\mathcal{K}(p^k,p^w,a)$ and $\mathcal{K}(p^k,p^w,a')$ are isomorphic if and only if
\begin{enumerate}
\item $a\equiv a'\; \pmod{p^w}$, for $k\geq 2w$; 
\item $a \equiv a'\; \pmod{p^{k-w}}$, for $2w> k> w$,
\item $a$ and $a'$ are any, for $w=k$.
\end{enumerate}
\end{cor}
\begin{proof}
Recall that, for the brace $B_{p^w}(p^k)$, the socle $\Soc(\Z_{p^k})=p^{k-w}\Z_{p^k}$ and note that the equation \eqref{eq:isosol} is equivalent to the following one
\[
a-a'\equiv a'p^wz-ap^{k-w}h \;\pmod{p^k},
\]
for some $h,z\in \Z_{p^k}$. 

Let $a-a'=p^wr$, for some $r\in\Z$. Then, by choosing $h=0$ and $z=r(a')^{-1}$, we obtain
\[a'p^wz-ap^{k-w}h=a'p^wr(a')^{-1}=p^wr=a-a'
\]
and $\mathcal{K}(p^k,p^w,a)\cong\mathcal{K}(p^k,p^w,a')$.

Let $a-a'=p^{k-w}r$, for some $r\in\Z$. Then, by choosing $h=-ra^{-1}$ and $z=0$, we obtain
\[a'p^wz-ap^{k-w}h=ap^{k-w}ra^{-1}=p^{k-w}r=a-a'
\]
and $\mathcal{K}(p^k,p^w,a)\cong\mathcal{K}(p^k,p^w,a')$ as well.

%Let $p^u|(a-a')$ for $u\in \{w,k-w\}$.  Then for some $r\in \Z$
%\[
%a-a'=p^ur=
%\begin{cases}
%a'p^wr(a')^{-1}-ap^{k-w}0\pmod{p^k}, & \text{for}\; u=w; \\
%a'p^w0-ap^{k-w}ra^{-1}\pmod{p^k}, & \text{for}\; u= k-w.
%\end{cases}
%\]
%This implies that .
 
Now assume that solutions $\mathcal{K}(p^k,p^w,a)$ and $\mathcal{K}(p^k,p^w,a')$  are isomorphic. Then
\[
a-a'\equiv a'p^wz-ap^{k-w}h \; \pmod{p^k}
\]
for some $h,z\in \Z_{p^k}$.

Thus if $k-w\geq w$ we have
$p^w|(a-a')$. Otherwise,  $p^{k-w}|(a-a')$. 
\end{proof}

It is known that finite cocyclic solutions are multipermutation solutions. In general it is difficult to describe the retraction of the solution. But in the case of indecomposable cocyclic solutions of the form $\mathcal{K}(p^k,p^w,a)$ we are able to give a precise description and also to compute their level.

\begin{theorem}\label{tw:ret}
Let $p$ be a prime number, $k,w\in \mathbb{N}$ and $a\in \Z_{p^k}^*$. For arbitrary $m\in \mathbb{N}$, the $m$-th retraction $\Ret^m(X)$ of the solution $(X,\sigma,\tau)=\mathcal{K}(p^k,p^w,a)$ is isomorphic to the solution:
 \[
 \mathcal{K}(p^{\max(k-mw,0)},p^w,a).
\]
\end{theorem}

\begin{proof}
%Let $(X,\sigma,\tau)=\mathcal{K}(p^k,p^w,a)$ and 
Suppose first $k\leq w$. Then $p^w\equiv 0\pmod {p^k}$ and the solution
$\mathcal{K}(p^k,p^w,a)$ is a permutation solution. Hence $\Ret^m(\mathcal{K}(p^k,p^w,a))$ is, for all $m>0$, the trivial solution~$\mathcal{K}(1,0,0)$.

Suppose now $k\geq w$.
At first we want to prove that $\mathrm{Ret}(X)\cong\mathcal{K}(p^{k-w},p^w,a)$. The relation~$\sim$ on~$X=\Z_{p^k}$ is defined as follows:
\[b\sim c \quad\Leftrightarrow \quad \sigma_b=\sigma_c\quad \Leftrightarrow \quad b\equiv c\pmod {p^{k-w}}.\]
Our planned isomorphism $\Phi:\mathrm{Ret}(X)\to \mathcal{K}(p^{k-w},p^w,a)$
simply sends $[b]_\sim$ to~$b \bmod {p^{k-w}}$. Such a mapping is clearly a bijection. Moreover, it is actually an identity on $\Z_{p^{k-w}}$
and both the operations $\sigma$ (in $\mathrm{Ret}(X)$ and in $\mathcal{K}(p^{k-w},p^w,a)$) are defined using the same $[a]_\sim$ and the same~$[p^w]_\sim$ and therefore $\Phi$ is a homomorphism.

The rest of the proof goes easily by an induction on $m$.
\end{proof}

\begin{cor}\label{prop:level}
The solution $\mathcal{K}(p^k,p^w,a)$ is a multipermutation solution of level~$\lceil k/w\rceil$.
\end{cor}

\begin{proof}
By Theorem~\ref{tw:ret}, the size of the retract $\Ret^m(X)$ is $\max(1,p^{k-mw})$.
And the smallest number~$m$ such that $k-mw\leq 0$ is $m=\lceil \frac kw\rceil$.
\end{proof}

\begin{ques}
We have presented two different solutions based on the brace $B_t(p^k)$, namely
an indecomposable solution $\mathcal{K}(p^k,t,a)$ and the (decomposable) associated solution. Nevertheless, according to Proposition \ref{prop:ring_1}, both the solutions have the same multipermutational level. We wonder whether this is a coincidence or there is a more general mechanism hidden behind.
\end{ques}

What remains is a description of indecomposable cocyclic solutions of size $4$. According to \cite{CPR20} or \cite{JPZ20} there are only two such solutions: one of multipermutational level~1 and one of multipermutational level~2. The latter one $(X,\tilde{\sigma},\tilde{\tau})$ can be constructed on the ring $(\Z_4,+,\cdot)$ using the operation
\[\tilde{\sigma}_b(c)=a+c+2ab,\]
for any chosen $a\in\{1,3\}$. We will denote such a solution by $\tilde{\mathcal{K}}(4,2,a)$. Similarly as in Remark \ref{rem:ess}, $\tilde{\mathcal{K}}(4,2,1)$ and $\tilde{\mathcal{K}}(4,2,3)$ are isomorphic via the mapping $\Phi(x)=3^{-1}x=3x$.

\subsection{Indecomposable cocyclic solutions of order $n=p_1^{k_1}\ldots p_s^{k_s}$}
To obtain finite  indecomposable solutions with the permutation group $\mathcal{G}(X)$ isomorphic, as a left brace, to a finite left brace $(B_t(n),\circ)$, for some pairwise distinct prime numbers $p_1,\ldots,p_s$, $k_1, \ldots, k_s\in \Z^+$ and $2<t\in \mathcal{N}(\Z_n)$, we can proceed exactly in the same way as in Subsection 3.1. 
There is one difference however---to get a trivial intersection of a family of subgroups of $(B_t(n),\circ)$ it is not necessary to assume that at least one of them is equal to $\{0\}$. But let us make such an assumption. In consequence, we obtain some indecomposable cocyclic solutions $(X,\sigma,\tau)$ defined on the set $X=\Z_{n}$ and being of the form  
$$\sigma_b(c):=a+c+tab+tac+t^2abc$$
for $b,c\in \Z_{n}$, $t\in \mathcal{N}(\Z_n)$ and $a\in \Z_{n}^*$. We will denote them by $\mathcal{K}(n,t,a)$.

On the other hand by \cite[Corollary 12]{CCS20} each indecomposable solution $(X,\sigma,\tau)$ of order $n=p_1^{k_1}\ldots p_s^{k_s}$ and with  cyclic permutation group $\mathcal{G}(X)$ is isomorphic to the product of indecomposable solutions $\mathcal{K}(p_i^{k_i},t_i,a_i)$, for $i=1,\ldots,s$.
A precise description of the isomorphism is given as follows:

\begin{theorem}\label{th:split}
Let $n=p_1^{k_1}\cdots p_s^{k_s}$ for some pairwise distinct prime numbers $p_1,\ldots,p_s$ and $k_1, \ldots, k_s\in \Z^+$. Let $t=t_1\cdots t_s\in\mathcal{N}(\Z_n)$, for $t_i\in \mathcal{N}(\Z_{p_i^{k_i}})$. Let $a\in \Z_n^*$ and let $a_i\in \Z_{p_i^{k_i}}^*$ be such that $a_i\equiv\frac{t}{t_i}a \pmod{p_i^{k_i}}$, for $1\leq i\leq s$. Then the solutions $\mathcal{K}(n,t,a)$ and 
$\mathcal{K}(p_1^{k_1},t_1,a_1)\times \mathcal{K}(p_2^{k_2},t_2,a_2)\times \ldots\times\mathcal{K}(p_s^{k_s},t_s,a_s)$ are isomorphic.
\end{theorem}
\begin{proof}
We show that the mapping $\Phi\colon \Z_n\to \Z_{p_1^{k_1}}\times \Z_{p_2^{k_2}}\times \ldots\times \Z_{p_s^{k_s}}$ 
given by
\[\quad x\mapsto (\frac{t}{t_1} x\pmod{p_1^{k_1}},\ldots,\frac{t}{t_i}x\pmod{p_i^{k_i}},\ldots,\frac{t}{t_s}x\pmod{p_s^{k_s}})
\]
is a homomorphism of the solutions. By Chinese Remainder Theorem it is also an isomorphism.

\begin{align*}
&\sigma_{\Phi(b)}\Phi(c)=\sigma_{(b\frac{t}{t_1}\pmod{p_1^{k_1}},\ldots ,b\frac{t}{t_s}\pmod{p_s^{k_s}})}
  (c\frac{t}{t_1}\pmod{p_1^{k_1}},\ldots ,c\frac{t}{t_s}\pmod{p_s^{k_s}})=\\
  &(\sigma_{b\frac{t}{t_1}\pmod{p_1^{k_1}}} (c\frac{t}{t_1}\pmod{p_1^{k_1}}),\ldots ,\sigma_{b\frac{t}{t_s}\pmod{p_s^{k_s}}} (c\frac{t}{t_s}\pmod{p_s^{k_s}}))=\\
&(a_1+\frac{t}{t_1}c+t_1\cdot a_1\cdot \frac{t}{t_1}b+t_1\cdot a_1\cdot\frac{t}{t_1}c+t_1^2\cdot a_1\cdot b\frac{t}{t_1}\cdot c\frac{t}{t_1}\pmod{p_1^{k_1}},\ldots )=\\
&(\frac{t}{t_1}a+\frac{t}{t_1}c+t\cdot a\cdot \frac{t}{t_1}b+t\cdot a\cdot\frac{t}{t_1}c+t_1\cdot a\cdot b\cdot c\frac{t^3}{t_1^2}\pmod{p_1^{k_1}},\ldots )=\\
&((a+c+tab+tac+t^2abc)\frac{t}{t_1}\pmod{p_1^{k_1}},\ldots,(a+c+tab+tac+t^2abc)\frac{t}{t_s}\pmod{p_s^{k_s}})=\\
&\Phi(a+c+tab+tac+t^2abc)=\Phi(\sigma_{b}(c)). \qedhere
\end{align*}
\end{proof}

\begin{corollary}\label{cor:wheniso}
Let $d:=gcd(t,\frac{n}{t})$. The solutions $\mathcal{K}(n,t,a)$ and $\mathcal{K}(n,t,a')$ are isomorphic if and only if
\[
a\equiv a'\; \pmod{d}.
\]
\end{corollary} 
\begin{proof}
We use here the notation of Theorem~\ref{th:split}. According to Corollary~\ref{cor:iso} and Theorem~\ref{th:split}, the solutions $\mathcal{K}(n,t,a)$ and $\mathcal{K}(n,t,a')$ are isomorphic
if and only if, for each $1\leq i\leq s$, the numbers $a_i=\frac {t}{t_i}a$ and $a_i'=\frac {t}{t_i}a'$, are congruent modulo $\min(t_i,p_i^{k_i}/t_i)$.
Since $\frac{t}{t_i}$ is coprime to $p_i$ and since $\min(t_i,p_i^{k_i}/t_i)=\gcd(t,\frac nt,p_i^{k_i})$, this is equivalent to
$a_i\equiv a_i'\pmod{\gcd(d,p_i^{k_i})}$, for each $1\leq i\leq s$.
According to Chinese Remainder Theorem, this means $a\equiv a'\pmod d$.
\end{proof}

Theorem \ref{th:split} describes all finite indecomposable cocyclic solutions except for some non-trivial ones of size $4n$, with $n$ odd. In this case we need to consider the four element solution which does not originate from a bicyclic brace. In the following proposition we give an explicit construction of the remaining solutions of size $4n$. We are not able to apply the map $\Phi\colon \Z_{4n}\to \Z_4\times \Z_n$ (from the proof of Theorem \ref{th:split}) directly since there is that four element cocyclic brace which is not cyclic. But we use a similar formula.

\begin{prop}\label{prop:4n}
Let $n$ be an odd number and $\bar 2$ be 
%the inverse of $2$ in the ring $(\Z_n,+,\cdot)$. 
a number such that $2\cdot \bar 2\equiv 1\pmod n$. 
Let $t\in \mathcal{N}(\Z_n)$ and $a\in \Z_{4n}^*$. Then $(\Z_{4n},\overline{\sigma},\overline{\tau})$ with 
\[ \overline{\sigma}_{b}(c)=a+c+2tab+4t\bar 2 ac+4t^2 abc\pmod{4n}\]
yields an indecomposable solution with its permutation group isomorphic to~$\Z_{4n}$.
\end{prop}

\begin{proof}
We will show that $(\Z_{4n},\overline{\sigma},\overline{\tau})$ is isomorphic to the product of two solutions: one defined on the set $\Z_4$ and the second one on $\Z_n$.

The mapping
\[\Phi\colon \Z_{4n}\to \Z_4\times \Z_n;\quad \Phi:x\mapsto (tx \pmod{4},2x \pmod{n})\]
gives for $b,c\in \Z_{4n}$
\begin{multline*}
\Phi(\overline{\sigma}_{b}(c))=(t(a+c+2tab+4t\bar 2 ac+4t^2 abc) \pmod{4},2(a+c+2tab+4t\bar 2 ac+4t^2 abc)\pmod{n})=\\
(ta+tc+2t^2ab \pmod{4},2a+2c+4tab+4tac+8t^2 abc\pmod{n})=\\
(ta+\Phi(c)+2ta\Phi(b) \pmod{4},2a+\Phi(c)+t2a\Phi(b)+t2a\Phi(c)+t^2 2a\Phi(b)\Phi(c)\pmod{n}).
\end{multline*}
Since for any $t\in \mathcal{N}(\Z_n)$ and $a\in \Z_{4n}^*$, $ta\in \{1,3\}\pmod{4}$ 
we obtain that $\Phi$ is a homomorphism between the solution $(\Z_{4n},\overline{\sigma},\overline{\tau})$ and a product of solutions $\tilde{\mathcal{K}}(4,2,ta\pmod{4})$ and
$\mathcal{K}(n,t,2a\pmod{n})$. Clearly, $\Phi$ is a bijection.
\end{proof}

\begin{theorem}\label{th:main}
A complete set of invariants for finite cocyclic indecomposable solutions are
\begin{itemize}
\item $n\in \N$;
\item $t\in \N$ such that 
\begin{itemize}
\item $t$ divides~$n$,
\item every prime~$p$ divides~$t$ whenever $p$ divides~$n$,
\item if $8$ divides~$n$ then $4$ divides~$t$;
\end{itemize}
\item $a\in\{1,\ldots,\mathrm{gcd}(t,n/t)\}$ coprime to~$n$.
\end{itemize}
\end{theorem}

\begin{proof}
Existence of such solutions is given in Theorem~\ref{th:split} and Proposition~\ref{prop:4n}. Moreover, we know that such a solution splits as a product of prime-power size solutions. 
If $n\equiv 4\pmod 8$ then there are two possibilities for the $2$-primary component and they depend on~$t$ only, not on~$a$.
Otherwise we use Corollary~\ref{cor:wheniso}.
\end{proof}

\begin{remark}\label{rem:ret}
Let $n$, $t$, $a$ be invariants described in Theorem \ref{th:main}. Similar proof as in Theorem \ref{tw:ret} shows that
$\Ret(\mathcal{K}(n,t,a))\cong \mathcal{K}(n/t,t,a).$
In particular, the solution $\mathcal{K}(n,t,a)$ is then a multipermutation solution of level $s$, where
$s$ is the smallest number such that $n$ divides $t^s$. Moreover, the solution $(\Z_{4n},\overline{\sigma},\overline{\tau})$ is of multipermutational level equal to $\max(2,s)$.
\end{remark}
\section{Enumeration and examples}

In the last section we compute the number of non-isomorphic indecomposable cocyclic solutions. 
We give exact numbers for prime powers only as for composed numbers the enumeration is straightforward using the spliting onto their prime components.
We, of course, omit the case $n=p$ where we have exactly one solution and the case $n=4$ where we have two solutions.

\begin{prop}\label{prop:howmany}
Let $p$ be a prime number and $k\in \Z^+$. Then the number $N$ of indecomposable solutions of order $p^k$ and with the cyclic  permutation group is equal to:
\begin{itemize}
\item $N=p^r+p^{r-1}-1$, 
for $p\neq 2$ and $k=2r\geq 2$
\item  $N=2p^r-1$, 
for $p\neq 2$ and $k=2r+1\geq 3$

\item $N=2^r+2^{r-1}-2$, 
for $p=2$ and $k=2r\geq 4$

\item $N=2^{r+1}-2$, 
for $p=2$ and $k=2r+1\geq 3$
\end{itemize}
\end{prop}
\begin{proof}
Directly by Corollary \ref{cor:iso}, for $p\neq 2$, we have
\begin{align*}
&N=\sum_{w=1}^{\lfloor \tfrac{k}{2}\rfloor}(p^w-p^{w-1})+\sum_{w=\lfloor \tfrac{k}{2}\rfloor+1}^{k-1}(p^{k-w}-p^{k-w-1})+1=\\
&p^{\left\lfloor \tfrac{k}{2}\right\rfloor}-1+
p^{k-\left\lfloor \tfrac{k}{2}\right\rfloor-1}-1+1=
p^{\left\lfloor \tfrac{k}{2}\right\rfloor}+
p^{k-\left\lfloor \tfrac{k}{2}\right\rfloor-1}-1.
\end{align*}
Moreover, for $p=2$
\begin{align*}
&N=\sum_{w=2}^{\lfloor \tfrac{k}{2}\rfloor}(2^w-2^{w-1})+\sum_{w=\lfloor \tfrac{k}{2}\rfloor+1}^{k-1}(2^{k-w}-2^{k-w-1})+1=\\
&2^{\left\lfloor \tfrac{k}{2}\right\rfloor}-2+
2^{k-\left\lfloor \tfrac{k}{2}\right\rfloor-1}-1+1=
2^{\left\lfloor \tfrac{k}{2}\right\rfloor}+
2^{k-\left\lfloor \tfrac{k}{2}\right\rfloor-1}-2.\qedhere
\end{align*}
\end{proof}

At the end we give some examples to illustrate Theorem~\ref{thm:constr}.

\begin{exm}\label{exm:n=9}
For $p=3$ and $k=2$ we have two non-isomorphic braces: $(B_3(9),+,\circ)$ and $(B_9(9),+,\circ)$. %It is evident that $\Z_9^*=\{1,2,4,5,7,8\}$.

Let us consider the brace $(B_3(9),+,\circ)$ with the multiplication 
$$x\circ y=x+y+3xy.$$
In this case $\Soc(B_3(9))=\{0,3,6\}$. Since $\Z^*_9=\{1,2,4,5,7,8\}$ then using the construction described in Theorem \ref{thm:constr} we obtain six indecomposable solutions:
\begin{align*}
&\mathcal{K}(9,3,1), \; {\rm with}\; \sigma_b(c)=1+3b+4c,\;
&\mathcal{K}(9,3,2), \; {\rm with}\; \sigma_b(c)=2+6b+7c, \\
& \mathcal{K}(9,3,4), \; {\rm with}\; \sigma_b(c)=4+3b+4c, \;
&\mathcal{K}(9,3,5), \; {\rm with}\;  \sigma_b(c)=5+6b+7c,\\
&\mathcal{K}(9,3,7), \; {\rm with}\; \sigma_b(c)=7+3b+4c, \;
&\mathcal{K}(9,3,8), \; {\rm with}\; \sigma_b(c)=8+6b+7c.
\end{align*}
By Corollary \ref{cor:iso}, the solutions $\mathcal{K}(9,3,1)$, $\mathcal{K}(9,3,4)$ and $\mathcal{K}(9,3,7)$ are isomorphic. In particular, the solutions $\mathcal{K}(9,3,1)$ and $\mathcal{K}(9,3,4)$ are isomorphic by the mapping $\Phi(x)=7^{-1}x=4x$, i.e. by the permutation $(147)(285)$.
Similarly, the solutions $\mathcal{K}(9,3,2)$, $\mathcal{K}(9,3,5)$ and $\mathcal{K}(9,3,8)$ are also isomorphic. 
Finally, again by Corollary \ref{cor:iso} solutions $\mathcal{K}(9,3,1)$ and $\mathcal{K}(9,3,2)$ are non-isomorphic. Note that they are isomorphic to solutions from Example \ref{exm:sol} for $r=3$ and $r=6$, respectively. The isomorphisms are given by permutations: $(285)$ and $(184572)(36)$, respectively.

For the trivial brace $(B_9(9),+,\circ)$ we have six isomorphic indecomposable solutions of multipermutational level~1: $\mathcal{K}(9,0,a)$, with $\sigma_b(c)=a+c$, for $a\in\{1,2,4,5,7,8\}$.
\end{exm}

By Corollary \ref{prop:level} and Proposition \ref{prop:howmany}
one can easily check how many indecomposable cocyclic solutions of each multipermutational level one can define on a set $X$ of cardinality $p^k$. Note that the highest possible level in this case does not exceed $k$. In particular, for each $2\neq p$ there exists at least one solution of level $k$, but for $p=2$ and $k>2$ there is none. Then the smallest example of the indecomposable solution of multipermutational level $3$ with cyclic permutation group can be constructed for $\left|X\right|=27$. 

\begin{exm}
There are five non-isomorphic indecomposable solutions with $\left|X\right|=27$ -- two of  multipermutational level $3$: $\mathcal{K}(27,3,1)$ and $\mathcal{K}(27,3,2)$ which originate from the brace $B_3(27)$, two of multipermutational level $2$: $\mathcal{K}(27,9,1)$ and $\mathcal{K}(27,9,2)$ which are obtained from the brace $B_9(27)$ and one of multipermutational level $1$: $\mathcal{K}(27,0,1)$ constructed from the trivial brace $B_{27}(27)$.
\end{exm}

Using Corollary \ref{cor:wheniso} and Remark \ref{rem:ret} we can quickly calculate the number of solutions and their permutational level in the case $n$ is product of powers of two primes.
\begin{remark}\label{rem:biprime}
Let $p$ and $q$ be two distinct prime numbers, $k,l\in \Z^+$ and $n=p^kq^l$. For $t=p^{k'}q^{l'}$ with $1\leq k'\leq k$, $1\leq l'\leq l$ there are $$(p^{\min(k',k-k')}-p^{\min(k',k-k')-1})(q^{\min(l',l-l')}-p^{\min(l',l-l')-1})$$ non-isomorphic solutions $\mathcal{K}(p^kq^l,p^{k'}q^{l'},a)$ with $a\in\Z_d^*$ and $d=p^{\min(k',k-k')}q^{\min(l',l-l')}$.

By Remark \ref{rem:ret} they are of multipermutational level equal to $\max(\lceil\dfrac{k}{k'} \rceil,\lceil\dfrac{l}{l'} \rceil)$.
\end{remark}

\begin{exm}
Consider a set $X$ with $|X|=3^3\cdot 5^2=675$. All indecomposable cocyclic solutions constructed on $X$ originate from bicyclic braces.  By Lemma \ref{lem:isobrace} there are 6 such braces (up to isomorphism). By Remark \ref{rem:biprime} there are 25 indecomposable cocyclic solutions and they are of multipermutational level at most~3. 

We focus now on multipermutation solutions of level~3 since there are not many examples in the literature. We have 10 non-isomorphic such solutions: 
$\mathcal{K}(675,3\cdot 5,a)$ for $a\in \{1,2,4,7,8,11,13,14\}$ and $\mathcal{K}(675,3\cdot 5^2,b)$ for $b\in \{1,2\}$.

In particular, by Theorem \ref{th:split} combined with Corollary \ref{cor:iso} we obtain the following sequence of isomorphisms:
\[
\mathcal{K}(27\cdot 25,3\cdot 5,1)\cong \mathcal{K}(27,3,5)\times \mathcal{K}(25,5,3)\cong\mathcal{K}(27,3,2)\times \mathcal{K}(25,5,3)\cong\mathcal{K}(27\cdot 25,3\cdot 5,76).
\]
By Remark \ref{rem:ret} we can describe now the retracts of $\mathcal{K}(27\cdot 25,3\cdot 5,1)$ in each step, i.e.\\ $\Ret(\mathcal{K}(27\cdot 25,3\cdot 5,1))\cong\mathcal{K}(45,15,1)$, $\Ret(\mathcal{K}(45,15,1))\cong\mathcal{K}(3,0,1)$ and $\left|\Ret(\mathcal{K}(3,0,1))\right|=1$.

Moreover we have $14$ non-isomorphic solutions of multipermutational level~$2$: $\mathcal{K}(675,3^2\cdot 5,a)$ for $a\in \{1,2,4,7,8,11,13,14\}$, $\mathcal{K}(675,3^2\cdot 5^2,b)$ for $b\in \{1,2\}$, $\mathcal{K}(675,3^3\cdot 5,c)$ for $c\in \{1,2,3,4\}$  and one solution $\mathcal{K}(675,0,1)$ of multipermutational level~1.
 \end{exm}

\begin{exm}
We construct all (up to isomorphism) indecomposable cocyclic solutions on set $X$ of cardinality $\left|X\right|=36=4\cdot 9$. By Remark \ref{rem:ret} they are multipermutation solutions of level at most~2. By \cite[Corollary 4.8]{JPZ20} there are $6$ such solutions. Using Theorem \ref{th:split} we obtain $3$ of them - those originating from bicyclic braces. There will be $1$ of multipermutational level~1: $\mathcal{K}(36,0,1)\cong \mathcal{K}(4,0,1)\times \mathcal{K}(9,0,4)\cong \mathcal{K}(4,0,1)\times \mathcal{K}(9,0,1)$. Note that $\mathcal{K}(36,0,1)\cong\mathcal{K}(36,0,a)$ for any $a\in \Z_{36}^*$. We also have $2$ non-isomorphic solutions of multipermutational level~2:
\begin{enumerate}
\item for $a\in\{1,7,13,19,25,31\}$: $\mathcal{K}(36,12,a)\cong \mathcal{K}(4,0,1)\times \mathcal{K}(9,3,1)$;
\item for $a\in\{5,11,17,23,29,35\}$: $\mathcal{K}(36,12,a)\cong \mathcal{K}(4,0,1)\times \mathcal{K}(9,3,2)$.
\end{enumerate} 
By Proposition \ref{prop:4n} we obtain $3$ non-isomorphic solutions with cyclic permutation group but not originating from the construction presented in Theorem \ref{thm:constr}. They are of multipermutational level~2.
\begin{enumerate}
\item for $t=9$ and any $a\in\Z_{36}^*$ we have a solution with $\overline{\sigma}_{b}(c)=a+c+18b \pmod{36}$. It is isomorphic to $\tilde{\mathcal{K}}(4,2,1)\times\mathcal{K}(9,0,1)$;
\item for $t=3$ and $a\in\{1,7,13,19,25,31\}$ we have a solution with $\overline{\sigma}_{b}(c)=a+25c+6b \pmod{36}$. It is isomorphic to $\tilde{\mathcal{K}}(4,2,1)\times\mathcal{K}(9,3,2)$;
\item for $t=3$ and $a\in\{5,11,17,23,29,35\}$ we have a solution with $\overline{\sigma}_{b}(c)=a+13c+30b \pmod{36}$. It is isomorphic to $\tilde{\mathcal{K}}(4,2,1)\times\mathcal{K}(9,3,1)$.
\end{enumerate}
\end{exm}

Proposition \ref{prop:howmany} together with Theorem \ref{th:split} and Proposition \ref{prop:4n} allow us to compute the precise number of all non-isomorphic indecomposable cocyclic solutions of arbitrary size $n$. Table \ref{Fig:count_all} shows that this number grows rapidly as $n$ increases. Comparing even with the number of all non-isomorphic left braces (see \cite[Table 5.4]{GV} for example for $n=25$ or $n=49$) one can see that these computations confirm Theorem \ref{th:main}.

\begin{table}[h!]
$$\begin{array}{|r|cccccccccccccc|}\hline
k & 2&3&4&5&6&7&8&9&10&11&12&13&14&15\\\hline
p=2&  2& 2& 4& 6& 10& 14 &22& 30& 46&62&94&126&190&  254\\
p=3&3 &  5& 11& 17& 35& 53& 107& 161& 323&485&971&1457&2915&4373\\
p=5& 5& 9& 29& 49& 149& 249& 749& 1249& 3749& 6249&18749&31249&93749&156249\\
p=7& 7& 13& 55& 97& 391& 685& 2743& 4801& 19207& 33613&134455&235297&941191&1647085\\\hline
\end{array}$$

\caption{The number of indecomposable cocyclic solutions of order $p^k$, up to isomorphism.}
\label{Fig:count_all}
\end{table} 

Finally, let us mention the infinite indecomposable solutions $(X,\sigma,\tau)$ with $\mathcal{G}(X)\cong \Z$. By  
\cite[Theorem 5.5]{CSV} a left  brace with multiplicative  group  isomorphic  to the group $\Z$ is trivial. Then, by \cite[Theorem 14]{CCS20} there is, up to isomorphism, unique cocyclic infinite indecomposable solution $(\Z,\sigma,\tau)$ of multipermutational level $1$ such that $\sigma_x(y)=y+1$. (See also \cite[Proposition 2]{JPZ20}.)


\begin{thebibliography}{99}

\bibitem{BCJ16} D. Bachiller, F. Cedó, E. Jespers, {\it Solutions of the Yang-Baxter equation associated with a left brace}, J. Algebra {\bf 463} (2016), 80--102.

\bibitem{CCP}
M. Castelli, F. Catino, G. Pinto, {\it Indecomposable involutive set-theoretic solutions of the Yang-Baxter equation}, J. Pure Appl. Algebra {\bf 223} (2019), 4477--4493.


\bibitem{CCS20} M. Castelli, F. Catino, P. Stefanelli, {\it Indecomposable involutive set-theoretic solutions of the
Yang-Baxter equation and orthogonal dynamical extensions of cycle sets},  available at {\tt http://arxiv.org/abs/2011.10083}

\bibitem{CPR20}
M. Castelli, G. Pinto, W. Rump, {\it On the indecomposable involutive set-theoretic solutions of the Yang-Baxter equation of prime-power size}, Comm. Algebra  {\bf 48} (2020), 1941--1955.

\bibitem{CGS17}
F. Ced\'{o}, T. Gateva-Ivanova, A. Smoktunowicz, {\it On the Yang–Baxter equation and left nilpotent left braces}, J. Pure Appl. Algebra {\bf 221} (2017), 751--756.

\bibitem{CGS18}
F. Ced\'{o}, T. Gateva-Ivanova, A. Smoktunowicz, {\it Braces and symmetric groups with special conditions}, J. Pure Appl. Algebra {\bf 222} (2018), 3877--3890.

\bibitem{CJO14}
F. Ced\'{o}, E. Jespers, J. Okni\'{n}ski, {\it Braces and the Yang-Baxter equation}, Comm. Math. Phys. {\bf 327} (2014), 101--116. Extended version arXiv:1205.3587.

%\bibitem{CJO}
%F. Ced\'{o}, E. Jespers, J. Okni\'{n}ski, {\it Primitive set-theoretic solutions of the Yang-Baxter equation}, available at {\tt http://arxiv.org/abs/2003.01983}


\bibitem{CSV}
F. Ced\'{o}, A. Smoktunowicz, L. Vendramin, {\it Skew left braces of nilpotent type}, Proc. Lond. Math. Soc. {\bf (3) 118(6)} (2019), 1367--1392.


\bibitem{Dr90}
V.G. Drinfeld, {\it On some unsolved problems in quantum group theory}, In: P.P. Kulish (ed.) Quantum
groups, in: Lecture Notes in Math., vol. 1510, Springer-Verlag, Berlin, (1992), pp. 1--8.

\bibitem{ESS}
P. Etingof,  T. Schedler, A. Soloviev, {\it Set-theoretical solutions to the quantum Yang-Baxter equation}, Duke Math. J. {\bf 100} (1999), 169--209.

\bibitem{GV}
L. Guarnieri, L. Vendramin, {\it Skew braces and the Yang-Baxter equation}, Math. Comp. {\bf 86} (307) (2017), 2519–-2534.


\bibitem{JPZ20}
P. Jedli\v cka, A. Pilitowska, A. Zamojska-Dzienio, {\it Indecomposable involutive solutions of the Yang-Baxter equation of multipermutational level 2 with abelian permutation group}, available at {\tt http://arxiv.org/abs/2011.00229}

\bibitem{Jimbo}
M. Jimbo, {\it Introduction to the Yang-Baxter equation}, Int. J. Modern Physics A, {\bf 4} (1989), 3759--3777.

\bibitem{KN07}
T. Kepka, P. N\v emec, {\it Commutative radical rings I}, Acta Univ. Carolinae Math. Phys. {\bf 48} (2007), 11--41.

\bibitem{KN08}
T. Kepka, P. N\v emec, {\it Commutative radical rings II}, Acta Univ. Carolinae Math. Phys. {\bf 49} (2008), 53--73.

\bibitem{Rump05}
W. Rump, {\it A decomposition theorem for square-free unitary solutions of the quantum
Yang-Baxter equation}, Adv. Math. {\bf 193} (2005), 40--55.

\bibitem{Rump07A}
W. Rump, {\it Braces, radical rings, and the quantum Yang-Baxter equation},
J. Algebra {\bf 307} (2007), 153--170.

\bibitem{Rump07B}
W. Rump, {\it Classification of cyclic braces}, J. Pure Appl. Algebra {\bf 209} (2007), 671--685.

\bibitem{Rump21}
W. Rump, {\it Cocyclic solutions to the Yang-Baxter equation}, Proc. Amer. Math. Soc. {\bf 149} (2021), 471--479.

\bibitem{Rump}
W. Rump, {\it Classification of indecomposable involutive set-theoretic solutions to the Yang–Baxter equation}, Forum Math. {\bf 32} (2020), 891--903.

\bibitem{S18}
A. Smoktunowicz, {\it On Engel groups, nilpotent groups, rings, braces and the Yang-Baxter equatiuon}, Trans. Am. Math. Soc. {\bf 370}(9) (2018), 6535--6564.

\bibitem{Vendr}
L. Vendramin, {\it Problems on skew left braces}, Adv. Group Theory Appl. {\bf 7} (2019), 15--37.
\end{thebibliography}
\end{document}